\documentclass{amsart}

\usepackage {amsfonts}
\usepackage[english]{babel}
\def\g{\gamma}
\def\G{\Gamma}
\def\d{\delta}
\def\a{\alpha}
\def\b{\beta}
\def\p{\varphi}
\def\e{\varepsilon}
\def\l{\lambda}
\def\L{\Lambda}
\def\o{\omega}
\def\O{\Omega}

\def\k{\kappa}

\def\D{\mathcal D}
\def\S{\mathcal S}
\def\cN{\mathcal N}

\def\R{{\mathbb R}}
\def\C{{\mathbb C}}
\def\N{{\mathbb N}}
\def\Z{{\mathbb Z}}

\DeclareMathOperator{\supp}{supp}
\DeclareMathOperator{\dist}{dist}

\newtheorem{Def}{Definition}
\newtheorem*{Th}{Theorem}
\newtheorem*{MTh}{Meyer's Theorem}

\begin{document}

\title{Application of methods of quasicrystals theory to entire functions of exponential growth}

\author{Sergii Yu.Favorov}

\address{Sergii Favorov,
\newline\hphantom{iii}  V.N.Karazin Kharkiv National University
\newline\hphantom{iii} Svobody sq., 4, Kharkiv, Ukraine 61022}
\email{sfavorov@gmail.com}

\maketitle {\small
\begin{quote}
\centerline {\bf Abstract.}

Let $f$ be an entire almost periodic function with zeros in a horizontal strip of finite width; for example, any exponential polynomial with purely imaginary exponents is such a function. Let $\mu$ be the measure on the set of zeros of $f$ whose masses coincide with multiplicities of zeros.  We define the Fourier transform in the sense of distributions for $\mu$ and prove that it is a pure point measure on $\R$ whose complex masses correspond to coefficients of Dirichlet series of the logarithmic derivative of $f$. Bases on this description and Meyer's theorem on quasicrystals, we give a simple necessary and sufficient condition for $f$ to be a finite product of sines.
\bigskip


\end{quote}
}

   \bigskip
\section{Introduction}\label{S1}
   \bigskip

In \cite{OU} A.Olevskii and A.Ulanovskii established  1-1 connection between zeros of real--rooted exponential polynomials and  Fourier quasicrystals. They proved that for the exponential polynomial
\begin{equation}\label{s}
P(z)=\sum_{1\le j\le N} q_j  e^{2\pi i\o_j z}, \qquad q_j\in\C,\quad \o_j\in\R,\quad N\in\N.
\end{equation}
with a zero set $A\subset\R$ the measure
\begin{equation}\label{a}
\mu_A=\sum_n \d_{a_n},\quad A=\{a_n\},
\end{equation}
is the Fourier quasicrystal, and conversely, every set $A\subset\R$ for which the measure \eqref{a} is a Fourier quasicrystal, is the zero set
of some exponential polynomial \eqref{s}. Note that previously a nontrivial example of the Fourier quasicrystal of the form \eqref{a} whose support has only a finite intersection with any arithmetical progression was found by P.Kurasov and P.Sarnak \cite{KS}.

Recall that a complex or real measure $\mu$ is called the {\it Fourier quasicrystal} if it is tempered distribution, has locally finite
support (that is, its intersection with any compact set is finite), its Fourier transform in the sense of distributions is also is a measure
with locally finite support, and moreover, both measures $|\mu|$, $|\hat\mu|$ are temperate distributions.
 Here and below $|\nu|(E)$ means the variation of the complex measure $\nu$ on the set $E$.

Recently, Fourier quasicrystals and connected with them real-rooted exponential polynomials are studied very actively (see, for example, \cite{M2,
LT, AKV}). In fact,  Fourier quasicrystals are the form of Poisson formulas (see below), the latter were used in particular by D. Radchenko and M. Viazovska in \cite{RV}.

L.Alon, A.Cohen, K.Vinzant \cite{ACV} established a connection between real--rooted exponential polynomials and Lee-Yang polynomials of many variables. In \cite{F3} the result of Olevskii and Ulanovskii was generalized to Dirichlet series, it was found a full description of zero sets of absolutely convergent real--rooted Dirichlet series
\begin{equation}\label{S}
Q(z)=\sum_{\o\in\O} q_\o  e^{2\pi i\o z}, \qquad q_\o\in\C,\quad \o\in\R,\quad\sum_{\o\in\O} |q_\o|<\infty,
\end{equation}
with bounded $\O$. F.Goncalves \cite{G} adapted the techniques of Olevskii and Ulanovskii for study of arbitrary tempered measures on $\R$
with a pure point Fourier transform and got similar results for such measures. In \cite{F5} we applied the above description of Fourier quasicrystals
 and Meyer's theorem on quasicrystals \cite{M} to prove a simple criterion for the exponential polynomials \eqref{s} and Dirichlet series \eqref{S} to be a finite product of sines. But there we assumed that these polynomials and series have only real zeros.

In \cite{F4} we  applied the method of Olevskii and Ulanovskii to  measures of the form \eqref{a} where  we replace the condition $A\subset\R$ with $A\subset\{z=x+iy:\,|y|\le H\}$. We slightly changed the definition of the Fourier transform of a measure, introduced
the measure $\hat\mu_A^c$, and  proved that the above results on zeros of exponential sums and Dirichlet series are valid in this case as well. Therefore it is natural to use this result to remove the condition of reality of zeros in \cite{F5}. Moreover, we would like to get the above criterion for the wider class of entire almost periodic functions.

To achieve this result, two obstacles must be overcome. First, the Dirichlet series of a holomorphic almost periodic function need not converge. Second, the inverse Fourier transform of the measure $\hat\mu_A^c$ is not the measure $\mu_A$, which makes direct application of Meyer's theorem impossible. Solving these problems, we obtain the following:

\begin{Th}
Let $f(z)$ be an entire  almost periodic function of exponential growth  with  zeros  in a strip
$$
S_{(\a,\b)}=\{z\in\C:\,\a<\Im z<\b\},\qquad -\infty<\a\le0\le\b<\infty,
$$
and
\begin{equation} \label{Q}
f'(z)/f(z)\sim\sum_{\g\in\G^*} h^+_\g e^{2\pi i\g z},\ z\in S_{(\b,\infty)},\qquad f'(z)/f(z)\sim\sum_{\g\in\G_*} h^-_\g e^{2\pi i\g z},\  z\in S_{(-\infty,\a)}.
\end{equation}
Then the condition
\begin{equation}\label{r2}
\sum_{|\g|<r} |h^+_\g|+\sum_{|\g|<r}|h^-_\g|<Kr\quad \forall\ r>1,
\end{equation}
where $K$ is a constant, is necessary and sufficient for the identity
\begin{equation}\label{sin}
f(z)=Ce^{iaz}\prod_{j=1}^J\sin(\a_j z+\b_j)
\end{equation}
 with some $C\in\C,\ J\in\N,\ a,\,\a_j,\,\b_j\in\R$.
\end{Th}
{\it Remark 1}. Obviously, this is  the case when support of $\mu$ is a finite union of arithmetic progressions.
There are many papers (e.g.,\cite{C}, \cite{F1}, \cite{LO}) where  various conditions  have been found  for support of $\mu$ to be just like that. But in these papers, the distances between points of  support of $\mu$ or between points of support of its Fourier transform $\hat\mu$  are assumed to be uniformly separated from zero.  The exception is Meyer's theorem \cite{M}, where support of $\mu$ is an arbitrary locally finite set and $\hat\mu$ is a slowly increasing Radon measure with arbitrary support. The proof of our theorem is based on it.

{\it Remark 2}. The representations \eqref{Q} for every entire almost periodic function of exponential growth with zeros in the strip
$S_{(\a,\b)}$ will be proved below.
\medskip

Almost all papers related to Fourier quasicrystals use (in explicit or implicit form) properties of almost periodic functions and measures.
In section \ref{S2} we give definitions of these objects and present some of their properties.
In section \ref{S3} we recall some definitions related to the Fourier transform of functions and measures.
In section \ref{S4} we give a proof of the theorem.

\bigskip
\section{Almost periodic functions and measures}\label{S2}
\bigskip

 Our proof of the Theorem is based on some results of the classical theory of almost periodic functions (see \cite{B}-\cite{B2}, \cite{Le}, \cite[Ch.6]{L}).

\begin{Def}\label{D1}
A continuous function $f(x)$ on the real line $\R$
is almost periodic  if for any
  $\e>0$ the set of $\e$-almost periods
  $$
E_{\e}= \{\tau\in\R:\,\sup_{x\in\R}|f(x+\tau)-f(x)|<\e\}
  $$
is relatively dense, i.e., $E_{\e}\cap(x,x+L)\neq\emptyset$ for all $x\in\R$ and some $L$ depending on $\e$;
a continuous function $g(z)$ on the open strip
$$
S_{(\a,\b)}=\{z=x+iy\in\C:\,-\infty\le\a<y<\b\le+\infty\}
$$
is almost periodic  if for any $\a',\,\b',\,\a<\a'<\b'<\b$, and
  $\e>0$ the set of $\e$-almost periods
  $$
E_{\a',\,\b',\e}= \{\tau\in\R:\,\sup_{x\in\R,\a'\le y\le\b'}|g(x+\tau+iy)-g(x+iy)|<\e\}
  $$
is relatively dense with $L$ depending on $\e,\a',\b'$.
\end{Def}

 Note that every exponential sum \eqref{s} and absolutely convergent Dirichlet series \eqref{S}
are almost periodic in the whole plane $S_{(-\infty,+\infty)}$.
\medskip

{\bf Properties of almost periodic functions}:
\smallskip

a) a finite sum, a product, and a linear combination of almost periodic  functions on $\R$ is also almost periodic on $\R$;
a finite sum, a product, and a linear combination of almost periodic  functions in a strip $S_{(\a,\b)}$ is also almost periodic in $S_{(\a,\b)}$,
\smallskip

b) each almost periodic function $f(x)$ on $\R$ is bounded and has unbounded support unless  $f(x)\not\equiv0$; each almost periodic function $g(z)$
on $S_{(\a,\b)}$ is bounded on every substrip $S_{(\a',\b')},\ \a<\a'<\b'<\b$, and has unbounded support unless  $g(z)\not\equiv0$,
\smallskip

c) the Fourier series corresponds to each  almost periodic function $f(x)$ on $\R$
$$
f(x)\sim\sum_{\o\in\O} a_\o e^{2\pi i\o x},
$$
where the spectrum $\O=\{\o:\,a_\o\neq0\}$  is a countable subset of $\R$;
the Fourier series corresponds to each  almost periodic function $g(z)$ in $S_{(\a,\b)}$
$$
g(x+iy)\sim\sum_{\o\in\O} a_\o(y) e^{2\pi i\o x}, \quad \a<y<\b,
$$
where the Fourier coefficients $a_\o(y)$ are continuous functions on $(\a,\b)$ and the spectrum $\O=\{\o:\,a_\o(y)\not\equiv0\}$  is a countable subset of $\R$; if, in addition, the function $g(z)$ is holomorphic on $S_{(\a,\b)}$, the Fourier series  turns into the Dirichlet series
$$
g(z)\sim\sum_{\o\in\O} h_\o e^{2\pi i\o z},
$$
where the coefficients $h_\o=a_\o(y) e^{-2\pi\o y}$ do not depend on $y$,
\smallskip

d) for every countable $\O\subset\R$  there are coefficients $\k_\o(n)$,
$$
0\le \k_\o(n)\le1,\qquad \k_\o(n)\to1\quad\text{as }\,n\to\infty,\qquad \#\{\o\in\O:\,\k_\o(n)\neq0\}<\infty\quad\text{for each fixed }\,n\in\N,
$$
such that for almost periodic $f(x)$ on $\R$ with spectrum in $\O$
$$
\sum_{\o\in\O} \k_\o(n) a_\o e^{2\pi i\o x}\to f(x)\quad (n\to\infty)
$$
uniformly with respect to $x\in\R$, and for almost periodic $g(z)$ on $S_{(\a,\b)}$ with spectrum in $\O$
\begin{equation}\label{BF}
\sum_{\o\in\O} \k_\o(n) a_\o(y) e^{2\pi i\o x}\to g(x+iy)\quad (n\to\infty)
\end{equation}
uniformly with respect to $x+iy\in S_{(\a',\b')},\ \a<\a'<\b'<\b,$ (Bochner--Feyer approximation).

\medskip
In the following items we assume that a function $g(z)$ is holomorphic in the  strips $S_{(\a,\b)},\ -\infty\le\a<\b\le\infty$.
\smallskip

e) If $g(z)$  almost periodic in a strip $S_{(\a,\b)}$ then its derivative $g'(z)$ is also almost periodic in the same strip, and its spectrum is a subset of spectrum $g$,
\smallskip

f) if $g(z)$ and $g_1(z)$ are  almost periodic in $S_{(\a,\b)}$ and their ratio $g(z)/g_1(z)$ is holomorphic in $S_{(\a,\b)}$, then it is almost periodic  in $S_{(\a,\b)}$,
\smallskip

g) if $g(z)$ is bounded and almost periodic in a strip $S_{(\a,\b)}$ with $\b=\infty$, then its  spectrum $\O$ is non-negative;
 if $g(z)$ is bounded and almost periodic  in a strip $S_{(\a,\b)}$ with $\a=-\infty$, then its spectrum $\O$ is non-positive,
\smallskip

h)if $g(z)$ is almost periodic in the strip $S_{(\a,\b)}$ with $\b=\infty$ and its spectrum $\O$ is non-negative, then the limit $\lim_{y\to\infty}g(x+iy)$ exists uniformly in $x\in\R$, and it is equal to the free term of the corresponding Dirichlet series; if $g(z)$ is almost periodic in the strip $S_{(\a,\b)}$ with $\a=-\infty$ and its spectrum $\O$ is non-positive,
then the limit $\lim_{y\to-\infty}g(x+iy)$ exists uniformly in $x\in\R$, and it is equal to the free term of the corresponding Dirichlet series,\smallskip

i) if $g(z)$ is an entire almost periodic function of exponential growth whose zeros are in some horizontal strip of finite width, then its spectrum $\O$ is bounded and $\inf\O\in\O,\ \sup\O\in\O$,
\smallskip

j) if $A$ is the zero set of an almost periodic in $S_{(\a,\b)}$ function $g(z)$ then for any $\e>0$
and $\a<\a'<\b'<\b$  there is $m=m(\e,\a',\b')>0$ such that the conditions $z\in S_{(\a'\b')}$ and $\dist(z,A)>\e$ implies
   $$
 |g(z)|\ge m.
 $$

k) for an almost periodic in $S_{(\a,\b)}$ function $g(z)$ the numbers
 $$
\#\{z=x+iy\in S_{(\a',\b')}\,:\,t\le x\le t+1,\,g(z)=0\},\quad \a<\a'<\b'<\b,
 $$
 are bounded uniformly in $t\in\R$; all zeros are counted according to their multiplicities, hence
 multiplicities of zeros  are uniformly bounded in every substrip $S_{(\a',\b')}$.
\medskip

Here we do not give a general definition of almost periodic measures on a strip (see \cite{R}, \cite{FRR}), but restrict ourselves to the case of measures on $\R$ (\cite{M1}).

\begin{Def}\label{D2}
A measure $\mu$ on $\R$ is almost periodic (in the sense of distributions) if the function
$$
\mu\star\phi(t)=\int\psi(t-x)\mu(dx)
$$
is almost periodic for any $\psi\in\D$.
\end{Def}
It is easily seen that any linear combinations of almost periodic measures is almost periodic, and the support of every almost periodic measure $\mu\not\equiv0$ is unbounded.

\bigskip
\section{Some properties of the Fourier transform}\label{S3}
\bigskip

Let $\D$ be the space of $C^\infty$-functions with compact support on $\R$, \  $\S$ be the Schwartz space of all $C^\infty$-functions $\p$ with finite norms
$$
\cN_k(\p)=\sup_{x\in\R}(1+|x|)^k\max_{0\le m\le k}|\p^{(m)}(x)|\quad \forall\, k=0,1,2,\dots,
$$
and the space of tempered distributions $\S'$ be the space  of all linear functional on $\S$, which continuous with respect to the topology generated by these norms.
The  Fourier transform
\begin{equation}\label{F}
   \hat\p(x)=\int_{\R}\p(t)e^{-2\pi i xt}dt
\end{equation}
is a continuous bijection of $\S$ onto $\S$, and for $\p\in\S$
\begin{equation}\label{IF}
 \p(t)=\int_{\R}\hat\p(x)e^{2\pi i xt}dx,
\end{equation}
(see, e.g., \cite{Ru}).

For every $\Phi\in\S'$ its Fourier transform $\hat\Phi$ is defined by the equality
\begin{equation}\label{h}
\hat\Phi(\p)=\Phi(\hat\p)\quad \forall\,\p\in\D.
\end{equation}
Since $\D$ is a dense subset of $\S$, we see that \eqref{h} is valid for all $\p\in\S$. Therefore, $\hat\Phi\in\S'$ for all $\Phi\in\S'$.
It follows from \eqref{F}--\eqref{h} that
$$
  \widehat{\hat\Phi}(x)=\Phi(-x).
$$
When $\Phi$ and $\hat\Phi$ are tempered measure, and
$$
\Phi=\sum_{\l\in\L}a_\l\d_\l,\qquad \hat\Phi=\sum_{\g\in\G}b_\g\d_\g
$$
 with countable $\L,\,\G$, then  equality \eqref{h} has the form
$$
   \sum_{\g\in\G}b_\g\p(\g)=\sum_{\l\in\L}a_\l\hat\p(\l),\qquad \forall\p\in\D,
$$
and is called {\it generalized Poisson formula}.

For $\p\in\D$  the function
\begin{equation}\label{ext}
  \hat\p^c(z)=\int_{\R}\p(t)e^{-2\pi i zt}dt=\widehat{(\p(t)e^{2\pi yt})}(x),\quad z=x+iy,
\end{equation}
 is the entire extension of $\hat\p$. It is easy to check that for every $k\in\N$ and $s>0$ uniformly in $y\in[-s,s]$
\begin{equation}\label{gro}
  \hat\p^c(x+iy)=O(|x|^{-k}),\qquad |x|\to\infty,
\end{equation}
In the case of $\supp\mu\subset\S_{(\a,\b)}, -\infty<\a<\b<\infty$, and $\mu$ such that
\begin{equation}\label{log}
\log\left(|\mu|\{z\in\S_{(\a,\b)}:\,|\Re z|<r\}\right)=O(\log r),
\end{equation}
 we defined the distribution $\hat\mu^c$ on $\R$ in \cite{F4} by the equality
$$
\hat\mu^c(\p)=\int_{\S_{(\a,\b)}}\hat\p^c(z)\mu(dz)\quad \forall\,\p\in\D.
$$

\bigskip
\section{The proof of the Theorem}\label{S4}
\bigskip

It is easy to check (see, e.g., \cite{F5}) that  functions of the form \eqref{sin} satisfy the condition \eqref{r2}. Prove the converse statement.

Let $f(z)$ be an entire almost periodic function of exponential growth with zeros in $\S_{(\a,\b)}$ for $-\infty<\a\le0\le\b<\infty$.
Since $f(z)$ has no zeros in the half-planes $\S_{(-\infty,\a)}$ and $\S_{(\b,\infty)}$,  then by the properties e) and f),
 the function $f'/f$ is almost periodic on these half-planes; by c), we get  representations \eqref{Q} there.

Let $\O$ be the spectrum of $f$.  By the property i), $\O$ is bounded and $\o_0=\inf_{\o\in\O}\o\in\O$. The function $f_1(z)=e^{-2\pi i\o_0 z}f(z)$ has non-negative spectrum, therefore by the property h), we have
$$
   \lim_{y\to\infty}f_1(x+iy)=h_{\o_0}\qquad\text{uniformly in}\ x\in\R.
$$
Since $\o_0\in\O$, we get that  $h_{\o_0}\neq0$. By e), the spectrum of $f'$ is a subset of spectrum $f$, hence the function $f_2(z)=e^{-2\pi i\o_0 z}f'(z)$ has non-negative spectrum too, and the limit
$$
   \lim_{y\to\infty}f_2(x+iy)
$$
exists and finite. Therefore the function $f'/f=f_2/f_1$ is bounded for $y$ large enough. It follows from the property g) that the set $\G^*$ from \eqref{Q} is a subset of $[0,+\infty)$. By the similar way, we get the set $\G_*$ is a subset of $(-\infty,0]$.

It follows from k) and j) that for sufficiently  small $\e$ and fixed $\eta>0$ there are $m=m(\e,\a,\b,\eta)>0$ and sequences $L_k\to+\infty,\,L'_k\to-\infty$  such that
$$
   |f(x+iy)|>m\quad\text{for}\quad x=L_k \quad\text{or}\quad x=L'_k,\quad \a-\eta\le y\le\b+\eta.
$$
Let $\p\in\D$, and $\hat\p^c(z)$ be defined in \eqref{ext}.
 Consider the integrals $I_k$ of the function $\hat\p^c(z)f'(z)f^{-1}(z)$ over the  boundaries of the rectangles
$$
\Pi_k=\{z=x+iy:\,L'_k<x<L_k,\,\a-\eta<y<\b+\eta\}.
$$
  Taking into account \eqref{gro}, we get that these integrals  tend to the difference
$$
   \int_{i(\a-\eta)-\infty}^{i(\a-\eta)+\infty}\hat\p^c(z)f'(z)f^{-1}(z)dz-
   \int_{i(\b+\eta)-\infty}^{i(\b+\eta)+\infty}\hat\p^c(z)f'(z)f^{-1}(z)dz=:I_*-I^*
$$
 as $k\to\infty$. On the other hand, $I_k$ tend to the sum
  \begin{equation}\label{r}
  2\pi i\sum_{\l:f(\l)=0}\text{Res}_\l \hat\p^c(z)f'(z)f^{-1}(z)=2\pi i\sum_{\l:f(\l)=0}a(\l)\hat\p^c(z)=2\pi i\int_{S_{(\a,\b)}}\hat\p^c(z)\mu_A(dz),
\end{equation}
where $a(\l)$ is the multiplicity of zero of $f(z)$ at the point $\l$, and $\mu_A=\sum_{\l:f(\l)=0}\d_\l$, where every zero $\l$ appears $a(\l)$ times. From k) it follows that the condition \eqref{log} is valid.

It follows from \eqref{Q}, \eqref{BF}, and \eqref{gro} that for some  $\k_\g^+(n),\,\k_\g^-(n)$ we get
\begin{multline*}
 I_*=\lim_{n\to\infty}\int_{-\infty}^{+\infty}\sum_{\g\in\G_*}\k_\g^-(n)h^-_\g e^{2\pi i\g(x+i(\a-\eta))}\hat\p^c(x+i(\a-\eta))dx\\
=\lim_{n\to\infty} \sum_{\g\in\G_*}\k_\g^-(n)h^-_\g e^{-2\pi\g (\a-\eta)}\int_{-\infty}^{+\infty}\hat\p^c(x+i(\a-\eta))e^{2\pi i\g x}dx,
\end{multline*}
\begin{multline*}
 I^*=\lim_{n\to\infty}\int_{-\infty}^{+\infty}\sum_{\g\in\G^*}\k_\g^+(n)h^+_\g e^{2\pi i\g(x+i\b+\eta)}\hat\p^c(x+i(\b+\eta))dx\\
=\lim_{n\to\infty} \sum_{\g\in\G^*}\k_\g^+(n)h^+_\g e^{-2\pi\g(\b+\eta)}\int_{-\infty}^{+\infty}\hat\p^c(x+i(\b+\eta))e^{2\pi i\g x}dx.
\end{multline*}
Combining \eqref{IF} and \eqref{ext}, we obtain
$$
  I_*=\lim_{n\to\infty}\sum_{\g\in\G_*\cap\text{supp }\p}\k_\g^-(n)h^-_\g\p(\g), \qquad I^*=\lim_{n\to\infty}\sum_{\g\in\G^*\cap\text{supp }\p}\k_\g^+(n)h^+_\g\p(\g).
$$
It follows from \eqref{r2} that coefficients $h_\g^\pm$ in \eqref{Q} satisfy the condition
$$
 \sum_{|\g|<r} |h^+_\g|+\sum_{|\g|<r}|h^-_\g|<\infty \quad\text{for all}\ r<\infty.
$$
Using Lebesgue's dominated convergence theorem and the property d), we get
$$
  I_*=\sum_{\g\in\G_*\cap\text{supp }\p}h^-_\g\p(\g), \qquad I^*=\sum_{\g\in\G^*\cap\text{supp }\p}h^+_\g\p(\g).
$$
Combining  \eqref{ext}, \eqref{h}, and \eqref{r}, we obtain
$$
\hat\mu_A^c(\p)=\int\hat\p^c(z)\mu_A(dz)=\frac{I^- -I^+}{2\pi i}=\sum_{\g\in\G^*\cap\text{supp }\p}\frac{ih^+_\g}{2\pi}\p(\g)-\sum_{\g\in\G_*\cap\text{supp }\p}\frac{ih^-_\g}{2\pi}\p(\g).
$$
By \eqref{r2}, we get for any   $r\in(1,\infty)$ and any $\p\in\D$ with support in $[-r,r]$
$$
  |\hat\mu^c_A(\p)|\le C\sup_\R |\p(t)|,\qquad C=C(r)<\infty.
$$
Therefore, the distribution $\hat\mu_A^c$ has an extension to a linear functional on the space of continuous functions $g$ on $[-r,r]$ such that $g(-r)=g(r)=0$ with the same estimate.
This means that $\hat\mu_A^c$ is a measure of the form
 \begin{equation}\label{ha}
\hat\mu_A^c=\sum_{\g\in\G^*}\frac{ih^+_\g}{2\pi}\d_\g-\sum_{\g\in\G_*}\frac{ih^-_\g}{2\pi}\d_\g.
\end{equation}

Set
 $$
 R(r)=\sum_{|\g|<r}|h_\g^+|+\sum_{|\g|<r}|h_\g^-|.
 $$
 By \eqref{r2}, we get
\begin{equation}\label{r1}
      R(r)=|\hat\mu_A^c|(-r,r)=O(r)\qquad (r\to\infty).
\end{equation}
Hence for any $\e>0$
$$
    \sum_{\g\in\G^*\setminus\{0\}}|h_\g^+|e^{-2\pi\g\e}\le \int_0^\infty e^{-2\pi r\e}R(dr)=2\pi \e\int_0^\infty R(r)e^{-2\pi r\e}dt<\infty.
$$
Therefore the Dirichlet series
$$
  \sum_{\g\in\G^*}h_\g^+e^{2\pi i\g z}
$$
absolutely converges on every half-plane $\Im z>\e$ for all $\e>0$. Hence this sum is the analytic extension of the function $f'(z)/f(z)$,
 and $f(z)\neq0$ for $y>0$.
By the same way we get $f(z)\neq0$ for $y<0$. So, $A\subset\R$, and the measure $\hat\mu_A^c$ coincides with the measure $\hat\mu_A$.

By \eqref{ha}, we also have for any $\p\in\D$
 \begin{equation}\label{c}
    \int\p(t-s)\mu_A(ds)=\int\hat\p(x)e^{2\pi itx}\hat\mu_A(dx)=\sum_{\g\in\G^*}\frac{ih^+_\g}{2\pi}\hat\p(\g)e^{2\pi it\g}-\sum_{\g\in\G_*}\frac{ih^-_\g}{2\pi}\hat\p(\g)e^{2\pi it\g}.
 \end{equation}
Since $\sup_{|\g|<r}|\hat\p(\g)|=O(r^{-2})$ as $r\to\infty$, we obtain that with some constant $C$.
$$
  \sum_{\g\in\G^*}|h^+_\g||\hat\p(\g)|+\sum_{\g\in\G_*}|h^-_\g||\hat\p(\g)|<C\left(R(1)+\int_1^\infty r^{-2}R(dr)\right)=2C\int_1^\infty R(r)r^{-3}dr<\infty.
$$
 By \eqref{r1}, the last integral is finite. Hence the both series in \eqref{c}  converge absolutely,  $\int\p(t-s)\mu_A(ds)$ is an almost periodic function, and $\mu_A$
 is an almost periodic measure. By k), masses of the measure $\mu_A$ take a finite number of integer values.
 
 The final part of the proof is very close to the same part in \cite{F5}. We provide it here for self-completeness.

Apply
\begin{MTh}{\rm(\cite[p.26]{M}, see also \cite{KL})}
If $\mu$ is a complex measure with locally finite support such that its masses take a finite number of values, and $\hat\mu$ is a measure such that
$$
|\hat\mu|(-r,r)=O(r),\qquad r\to\infty,
$$
then each set $\L_s=\{\l\in\R:\,\mu(\l)=s\}$ belongs to the coset ring\footnote{The coset ring of an abelian group $G$ is the smallest collection of subsets of $G$ which is closed under finite unions,
finite intersections, and complements  and which contains all shifts of all subgroups of $G$.} of $\R$.
\end{MTh}
Let $K=\max\{k:\,\L_k\neq\emptyset\}$. Then each set $\L_k=\{\l:\,a_\l=k\},\,1\le k\le K$, belongs to the coset ring of $\R$. H.Rosenthal \cite{Ro} proved
that every discrete (without finite limit points) set $\L\subset\R$  belonging to the coset ring is a finite union of arithmetic progressions up to a finite set. So
\begin{equation}\label{p}
  \L_k=\left[\bigcup_{m=1}^{M_k} L_{k,m}\cup F_k^+\right]  \setminus F_k^-,\quad L_{k,m}=\{\a_{k,m}n+\b_{k,m}:\,n\in\Z\},\quad 1\le k\le K,\,\,1\le m\le M_k,
\end{equation}
where $F_k^+,\, F_k^-$ are finite. If the intersection of two progressions $L=\{\a n+\b:\,n\in\Z\}$ and $L'=\{\a'n+\b':\,n\in\Z\}$ contains more than one point, then $\a/\a'$ is rational and all  sets
$L\setminus L',\,L'\setminus L,\,L\cap L',\,L\cup L'$ are finite unions
of disjoint arithmetic progressions of the type $\{\a\a'n+\b'':\,n\in\Z\}$. Therefore we can transform \eqref{p} such  that  any distinct $L_{k,m}$ and $L_{k'm'}$  have at most one common point.
Set
$$
\mu_{k,m}=\sum_{\l\in L_{k,m}} k\,\d_\l,\quad 1\le k\le K,\quad 1\le m\le M_k.
$$
Since these measures are periodic, we get that the measure
$$
  \nu=\mu_A-\sum_{k=1}^K\sum_{m=1}^{M_k} \mu_{k,m}
$$
is almost periodic. Its support is contained in the finite set
$$
\bigcup F_k^+\cup F_k^-\cup(L_{k,m}\cap L_{k'm'}),
$$
where the union is taken over all $1\le k,k'\le K,\,1\le m,m'\le M_k$ provided that pairs $(k,m)$ and $(k',m')$ do not coincide. Therefore,
$\nu\equiv0$ and
$$
\mu_A=\sum_{k=1}^K\sum_{m=1}^{M_k} \mu_{k,m}.
$$
Hence the function

$$
 D(z)= \frac{f(z)}{\prod_{k=1}^K\prod_{m=1}^{M_k}\sin^{k}(\pi z/\a_{k,m}-\pi\b_{k,m}/\a_{k,m})}
$$
\medskip

\noindent is entire without zeros. By \cite[Ch.1]{L},  $D(z)$ is a function of exponential growth and is equal to $Ce^{iaz}$. By property f), this function is  almost periodic in the whole plane $S_{(-\infty,\infty)}$, hence it is bounded on the line $l=\{z=t+i:\,t\in\R\}$. Consequently, $a\in\R$ and we obtain \eqref{sin}. The Theorem is proved.

\end{document}